\documentclass[11pt]{article}
\usepackage{amsmath, amsthm, amssymb}
\begin{document}
\title{The Weil-Petersson Geometry On the Thick Part of the Moduli 
Space of Riemann Surfaces}

\author{Zheng Huang}
\date{Feb. 17, 2006}
\newtheorem{theorem}{Theorem}[section]

\newtheorem{pro}[theorem]{Proposition}
\newtheorem{cor}[theorem]{Corollary}
\newtheorem{lem}[theorem]{Lemma}
\newtheorem{rem}[theorem]{Remark}

\newcommand{\WP}{Weil-Petersson}
\newcommand{\TS}{Teichm\"{u}ller space}
\newcommand{\ms}{moduli space}
\newcommand{\cs}{conformal structure}
\newcommand{\hqd}{holomorphic quadratic differential}
\newcommand{\RS}{Riemann surface}
\newcommand{\hm}{harmonic map}
\newcommand{\hym}{hyperbolic metric}
\newcommand{\Sc}{sectional curvature}
\newcommand{\Rc}{Ricci curvature}
\newcommand{\Jd}{Jacobian determinant}
\newcommand{\cd}{complex dimension}
\newcommand{\Bd}{Beltrami differential}
\newcommand{\ts}{tangent space}
\newcommand{\mcg}{mapping class group}
\newcommand{\af}{asymptotic flatness} 
\newcommand{\hsc}{holomorphic sectional curvature}

\maketitle
\footnotetext{Zheng Huang: Department of Mathematics, 
University of Michigan, Ann Arbor, MI 48109, USA. Email: zhengh@umich.edu}
\tableofcontents
\section {Introduction}

Results on asymptotic geometry of the {\WP} metric on {\TS} are presented 
in this paper. However, there is \emph{no curve pinching} in the setting of 
this paper. We restrict ourselves on the thick part of the {\ms}, where 
there is a positive lower bound on the least hyperbolic length of any short 
closed curves on the surface. We obtain boundedness of the 
Riemannian curvatures on the thick part of the {\ms} with respect to the 
genus $g$, by considering {\bf surfaces with large genera}.

The motivation of this paper is the recent consideration of articles 
\cite{Mi1} and \cite{Mi2} of Maryam Mirzakhani, where she studied the 
hyperbolic lengths of simple closed geodesics on a {\RS} and calculated 
{\WP} volume of the {\ms} of {\RS}s. The {\WP} volume grows exponentially 
in $g$ (see also {\cite{Wp2}}), in particular, a rational multiple 
of ${\pi}^{6g-6+2n}$, where $n$ is the number of punctures on the surface. 
Thus it is natural to ask how the curvatures vary in terms of $g$. 
We find, in theorem 1.2, that the absolute value of the {\WP} {\Sc} is bounded 
with respect to $g$, in the thick part of the {\ms}.

Let $\Sigma$ be a smooth, oriented, closed {\RS} of genus $g$, and we always 
assume $g > 1$. The study of {\cs}s, or equivalently {\hym}s, on a compact 
{\RS} naturally leads to the classical Teichm\"{u}ller theory.

{\TS} ${\mathcal {T}}_g$ is the space of {\cs}s on $\Sigma$, modulo an 
equivalent relationship, where two {\cs}s $\sigma$ and $\rho$ are equivalent 
if there is a biholomorphic map between $(\Sigma,\sigma)$ and $(\Sigma,\rho)$, 
in the homotopy class of the identity. Riemann's {\ms} $\mathcal{M}_g$ of 
{\RS}s is obtained as the quotient of {\TS} by the {\mcg}.

{\TS} ${\mathcal {T}}_g$ is known (\cite {Ah1}) to be a complex manifold 
of {\cd} $3g - 3$,  and the co{\ts} at $\Sigma$ is identified with 
$QD(\Sigma)$, the space of {\hqd}s. One also identifies the {\ts} 
with $HB(\Sigma)$, the space of harmonic {\Bd}s.

We recall that the {\WP} cometric on {\TS} is defined on $QD(\Sigma)$ by 
the $L^2$-norm:
\begin{eqnarray}
||\phi||_{WP}^2 = \int_{\Sigma} \frac {|\phi|^2}{\sigma}dzd\bar{z},
\end{eqnarray}
where $\sigma |dz|^2$ is the {\hym} on $\Sigma$. By duality, we obtain a 
Riemannian metric on the {\ts} of ${\mathcal{T}}_g$: for $\mu \in HB(\Sigma)$, 
we have 
\begin{eqnarray}
||\mu||_{WP}^2 = \int_{\Sigma}|\mu|^2 {\sigma}dzd\bar{z} = \int_{\Sigma}|\mu|^2 dA. 
\end{eqnarray}

The differential geometry of this metric on {\TS} has been extensively studied. 
Ahlfors showed that the {\Rc}, the {\hsc} and scalar curvatures are negative 
(\cite {Ah1}, \cite {Ah2}). Royden showed the {\hsc} is bounded away from 
zero, and conjectured an explicit bound ${\frac{-1}{2 \pi (g-1)}}$ (\cite {R}), 
later proved by Wolpert (\cite {Wp3}).

Throughout this paper, we denote $l_{0}(\sigma)$ as the {\bf systole}, the 
length of the shortest closed geodesic, of the surface $\Sigma$. As an important 
invariant of {\hym}s, the systole is related to several important aspects of the 
{\WP} geometry, such as the injectivity radius and curvatures near the 
compactification divisor, as seen in {\cite {Wp4}, \cite {Wp5}, \cite {H1}, \cite {H2}}.

In this paper, we assume there is a positive lower bound $r_0$ on the injectivity radius, 
denoted by $inj_{\sigma}(\Sigma)$, of $\Sigma$, therefore the systole $l_{0}(\sigma)$ also 
has a positive lower bound. We call {\bf{the thick part}} of the {\ms} as the compact subset 
of the {\ms} where $inj_{\sigma}(\Sigma) > r_0 > 0$, i.e., {\emph{the thick part of {\ms} 
consists of fat surfaces}}.

Our first result is concerning the {\WP} {\hsc} on a large part of the {\ms}:
\begin{theorem}
There exists a positive constant $C_1$, independent of $g$, such that the {\WP}
{\hsc} $K_h$ satisfies that 
\begin{center}
$-{C_1} < K_h < {\frac{-1}{2 \pi (g-1)}}$,
\end{center}
on the thick part of the {\ms}.
\end{theorem}

Tromba and Wolpert proved that the {\WP} {\Sc}s are negative (\cite {Tr}, \cite 
{Wp3}). For the asymptotics of the {\Sc}s, much of the consideration has been on 
the {\WP} geometry near the compactification divisor. It is shown that, there are 
directions (near the compactification divisor) in which the {\Sc}s have no negative 
upper bound (\cite {H1}), nor lower bound (\cite {H2}, \cite {Tp}), for fixed $g$.

When the surface is away from the compactification divisor, from a compactness 
argument, all curvatures are bounded and the bounds are in terms of $g$ and 
$l_{0}(\sigma)$. In this paper, we find the curvatures are bounded, independent of 
the genus $g$, i.e.,
\begin{theorem}
On the thick part of the {\ms}, there is a positive constant $C_2$, independent of the 
genus $g$, such that the {\Sc} $K$ of the {\WP} metric satisfies that 
\begin{center}
$-{C_2} < K < 0$.
\end{center}
\end{theorem}

As an easy application, we also find that the following estimates of the {\Rc} 
and the scalar curvature:
\begin{theorem}
On the thick part of the {\ms}, there are positive constants $C_3, C_{4}$, 
independent of the genus $g$, such that the {\Rc} of the {\WP} metric satisfies that 
\begin{center}
$-C_3 g< Ric(d_{WP}) < {\frac{-1}{2 \pi (g-1)}}$.
\end{center}
and the scalar curvature $s$ satisfies that 
\begin{center}
$-C_4 g^2 < s < -{\frac{9g-6}{4\pi}}$.
\end{center}
\end{theorem}

The main technique involved is the analysis of {\hm}s between hyperbolic 
surfaces. The influence of Eells-Sampson's pioneer work (\cite {ES}) is reflected in 
many perspectives of Teichm\"{u}ller theory, for example, see \cite {J}, 
\cite {EE}, \cite {Wf1}, \cite {Wf2}, \cite{Min}, and many others. One 
often finds that many functions associated to {\hm}s satisfy some 
geometric differential equations, where one can apply techniques in geometric analysis.

Here is the outline of the paper. In $\S 2$, we give a quick exposition of the 
background and introduce our notions. We prove our main theorems in $\S 3$.

The author is grateful to Maryam Mirzakhani for bringing this problem to his 
attention, and Xiaodong Wang, Mike Wolf, for their very insightful comments 
and generous help.

\section{Background}

In this paper, let $\Sigma$ be a fixed, compact, smooth, oriented surface of genus 
$g > 1$. We also denote $\sigma |dz|^2$ as a {\hym} on $\Sigma$, for conformal 
coordinates $z$. On $(\Sigma, \sigma |dz|^2)$, we denote the Laplacian as
\begin{center}
$\Delta = {\frac{4}{\sigma}}{\frac{\partial^2}{\partial z \partial \bar{z}}}$,
\end{center}
with nonpositive eigenvalues.

{\TS} is a complex manifold with the {\mcg} acting by biholomorphisms, of {\cd} 
$3g-3$, the number of independent closed curves in any pair-of-pants decomposition 
of the surface.

The {\WP} metric is invariant under the action of the {\mcg}, hence it decends to a 
metric on the {\ms}. It is shown that the {\WP} metric is K\"{a}hlerian (\cite {Ah1}), 
with negative {\Sc} (\cite {Tr}, \cite {Wp3}). The {\WP} Riemannian curvature 
tensor is given as the Tromba-Wolpert formula
(\cite {Tr}, \cite {Wp3}):
\begin{center} 
$R_{\alpha\bar{\beta}\gamma\bar{\delta}} = 
\int_{\Sigma}D(\mu_{\alpha}\bar\mu_{\beta})\mu_{\gamma}\bar\mu_{\delta}dA + 
\int_{\Sigma}D(\mu_{\alpha}\bar\mu_{\delta})\mu_{\gamma}\bar\mu_{\beta}dA$,
\end{center}
where $\Delta$ is the Laplacian for the {\hym} $\sigma$ on the surface, and 
$D = -2(\Delta -2)^{-1}$ is a self-adjoint, compact operator.
Those $\mu$'s in the formula are tangent vectors, i.e., harmonic {\Bd}s. There is a natural 
pairing between 
$QD(\Sigma) =\{\phi(z) dz^2\}$ and $HB(\Sigma) =\{\mu(z) {\frac{d\bar{z}}{dz}}\}$:
\begin{center}
$<\phi dz^2,\mu {\frac{d\bar{z}}{dz}}> = Re \int_{\Sigma}\phi \mu dzd\bar{z}$.
\end{center}

A remarkable property of the {\WP} metric is that it is incomplete (\cite {Ma}, 
\cite {Wp1}). This is caused by pinching off at least one short closed geodesic on the 
surface. The {\WP} completion modulo the {\mcg} is topologically the Deligne-Mumford 
compactification of the {\ms}. The compactification divisor, thus 
consists of a union of lower dimensional {\TS}s, each such space consists of noded 
{\RS}s, obtained by pinching nontrivial short closed geodesics on the surface 
(\cite {B}, \cite {Ma}). Therefore the compactification divisor can be described 
via the systole $l_{0}(\sigma)$ as the set $\{l_{0}(\sigma) = 0\}$.

Since we will analyze {\hm}s between compact hyperbolic surfaces, we recall some 
fundamental facts here. For a Lipschitz map $w:(\Sigma, \sigma |dz|^2) \rightarrow 
(\Sigma, \rho |dw|^2)$, where $\sigma |dz|^2$ and $\rho |dw|^2$ are {\hym}s on 
$\Sigma$, and $z$ and $w$ are conformal coordinates on $\Sigma$, 
one follows Sampson (\cite {S}) to define 
\begin{center}
${\mathcal{H}}(z) = {\frac{\rho(w(z))}{\sigma(z)}}|w_z|^2, {\mathcal{L}}(z) = 
{\frac{\rho(w(z))}{\sigma(z)}}|w_{\bar{z}}|^2$.
\end{center}
We call ${\mathcal{H}}(z)$ the {\it holomorphic energy density}, and ${\mathcal{L}}(z)$ 
the {\it anti-holomorphic energy density}. Then the energy density function of $w$ is 
simply $e(w)= {\mathcal{H}} + {\mathcal{L}}$, and the total energy is then given by
\begin{center}
$E(w,\sigma,\rho) = \int_{\Sigma}e\sigma|dz|^2$.
\end{center}

We also note that the {\Jd} relative to the $\sigma$ metric is therefore given by 
$J(z) = {\mathcal{H}}(z) - {\mathcal{L}}(z)$.

The map $w$ is called {\it harmonic} if it is a critical point of this energy 
functional, i.e., it satisfies Euler-Lagrange equation:
\begin{center}
$w_{z\bar{z}}+ {\frac{\rho_w}{\rho}}w_z w_{z\bar{z}} = 0$.
\end{center}

The $(2,0)$ part of the pullback $w^{*}\rho$ is the so-called {\it {Hopf differential}}:
\begin{center}
$\phi(z)dz^2 = (w^{*}\rho)^{(2,0)} = \rho(w(z)) w_z {\bar{w}}_zdz^2$.
\end{center}
It is routine to check that $w$ is harmonic if and only if $\phi dz^2 \in QD(\Sigma)$, 
and $w$ is conformal if and only if $\phi=0$.

In our situation, there is a unique 
{\hm} $w:(\Sigma, \sigma) \rightarrow (\Sigma, \rho)$ in the homotopy class of the 
identity, moreover, this map $w$ is a diffeoemorphism with positive Jacobian $J$, and 
${\mathcal{H}}>0$ (\cite {ES}, \cite {Hr}, \cite {S}, \cite {SY}).

A key observation to link the harmonic maps to Teichm\"{u}ller theory is that one 
obtains a map from {\TS} to $QD(\Sigma)$, for some fixed {\hym} $\sigma$. More 
specifically, this map sends any {\hym} on $\Sigma$ to a holomorphic quadratic 
differential associated to the unique {\hm} in the homotopy 
class of the identity. This map is a diffeomorphism (\cite {S}, \cite{Wf1}).

\section{Proof of Main Theorems} 
\noindent 
Let $\mu = \mu (z){\frac{d\bar{z}}{dz}}$ be a unit {\WP} normed harmonic {\Bd}, thus 
$\int_{\Sigma}|\mu|^2 dA = 1$. The {\WP} {\hsc} in the direction of $\mu$ is given by
\begin{center}
$K_h = -2\int_{\Sigma}D(|\mu|^2)|\mu|^2 dA$.
\end{center} 
Its upper bound, in terms of the genus, is known, proved by Wolpert:
\begin{lem}(\cite {Wp4})
$K_h < -{\frac{1}{2\pi(g-1)}}$.
\end{lem}

\begin{rem}
This upper bound holds without restriction on the systole. Since all 
{\Sc}s are negative, so the {\Rc} and scalar curvature are bounded from 
above by $-{\frac{1}{2\pi(g-1)}}$, and $-{\frac{3(3g-2)}{4\pi}}$, respectively (\cite {Wp4}).  
\end{rem}

We now show the following pointwise estimate on $|\mu(z)|$, where the tangent 
vector $\mu(z){\frac{d\bar{z}}{dz}}$ is normalized to have unit {\WP} 
norm, and we shall apply this estimate to prove our main theorems.
\begin{theorem}
For $\mu(z){\frac{d\bar{z}}{dz}} \in HB(\Sigma)$ with $||\mu||_{WP} =1$, there 
exists a positive constant $h_0$, independent of $g$, such that 
$|\mu (z)| \le h_0$, for all $z \in \Sigma$, where the surface $\Sigma$ is in 
the thick part of the {\ms}.
\end{theorem}
\begin{proof}
Recall that $\mu(z){\frac{d\bar{z}}{dz}}$ is a harmonic {\Bd}, hence is a symmetric 
tensor given as $\bar{\phi}(ds^2)^{-1}$ for $\phi$ a {\hqd} with at most simpole poles 
at the cusps and $ds^2$ the {\hym} tensor (\cite{Wp2}). Since the surface $\Sigma$ lies in the 
thick part of {\ms}, hence no cusps and this {\hqd} $\phi = \phi(z)dz^{2}$ has no poles.

Note that $\phi \in QD(\Sigma)$, as stated in the previous section, by a theorem of 
Wolf (\cite{Wf1}), there exists a {\hym} $\rho$ on surface $\Sigma$, and a unique 
{\hm} $w: (\Sigma,\sigma) \rightarrow (\Sigma,\rho)$, such that $\phi$ is the 
Hopf differential associated to this {\hm} $w$, i.e., 
$\phi(z)dz^{2} = \rho(w(z)) w_z {\bar{w}}_zdz^{2}$.

Much of our study will be analyzing this {\hm} $w$. Note that even though 
$inj_{\sigma}(\Sigma) > r_{0} >0$, the metric $\rho$ might not lie in the thick part 
of the {\ms}. We recall that the holomorphic and anti-holomorphic energy density 
functions of $w$ are defined as 
${\mathcal{H}}(z) = {\frac{\rho(w(z))}{\sigma(z)}}|w_z|^2$, and ${\mathcal{L}}(z) = 
{\frac{\rho(w(z))}{\sigma(z)}}|w_{\bar{z}}|^2$, respectively.

The energy density function of $w$ is $e(w)= {\mathcal{H}} + {\mathcal{L}}$, while the {\Jd} 
between {\hym}s $\sigma$ and $\rho$ is therefore 
$J(z) = {\mathcal{H}}(z) - {\mathcal{L}}(z)$. Since the map $w$ is a diffeomorphism 
with positive {\Jd}, we have ${\mathcal{H}}(z) > {\mathcal{L}}(z) \ge 0$.

We also find that 
${\mathcal{H}}{\mathcal{L}} = {\frac{|\phi|^{2}}{\sigma^{2}}} = |\mu|^{2}$, so 
the zeros of ${\mathcal{L}}$ are the zeros of $|\mu|$, or equivalently, the zeros of $\phi$.

Let $\nu$ be the {\Bd} of the map $w$, defined by 
$\nu = {\frac{w_{\bar z}d\bar z}{w_{z}dz}}$. It measures the failure of $w$ to be 
conformal, and since $J > 0$, we have $|\nu|<1$.

One easily finds that
 $|\nu|^{2} = {\frac{\mathcal{L}} {\mathcal{H}}}$, therefore,
 \begin{center}
 $|\mu|= {\sqrt{{\mathcal{H}}{\mathcal{L}}}} = {\mathcal{H}}|\nu| < {\mathcal{H}}$.
 \end{center}
Thus it suffices to estimate ${\mathcal{H}}(z)$ to bound $|\mu|$ pointwisely.

Let $z_{0} \in \Sigma$ such that ${\mathcal{H}}(z_{0}) = max_{z \in \Sigma}{\mathcal{H}}(z)$. 
We follow a calculation of Schoen-Yau to define a local one-form $\theta = \sqrt{\sigma(z)}dz$, 
and find (\cite {SY}):
 \begin{center}
 $|w_{\theta}|^{2}= {\frac{\rho}{\sigma}}|w_{z}|^{2} = {\mathcal{H}}(z)$,
 \end{center}
 and 
 \begin{eqnarray}
 \Delta |w_{\theta}|^{2}= 4|w_{\theta \theta}|^{2}+2J|w_{\theta}|^{2}-2|w_{\theta}|^{2}.
 \end{eqnarray}
 We rewrite this as 
 \begin{eqnarray}
 \Delta{\mathcal{H}}= 4|w_{\theta \theta}|^{2}+2J{\mathcal{H}}-2{\mathcal{H}} > -2{\mathcal{H}}.
 \end{eqnarray}

 Therefore ${\mathcal{H}}$ is a subsolution to an elliptic equation $(\Delta+2)f=0$.

 Recalling that $inj_{\sigma}(\Sigma)>r_{0}>0$, we embed a hyperbolic ball $B_{z_{0}}({\frac{r_{0}}{2}})$ 
 into $\Sigma$, centered at $z_{0}$ with radius ${\frac{r_{0}}{2}}$. Morrey's theorem (\cite{Mo}, theorem 5.3.1) 
 on subsolutions of elliptic differential equations guarantees that there is a constant $C(r_{0})$, such that, 
 \begin{center}
 ${\mathcal{H}}(z_{0}) = sup_{B_{z_{0}}({\frac{r_{0}}{4}})}{\mathcal{H}}(z) 
 \le C(r_{0})\int_{B_{z_{0}}({\frac{r_{0}}{2}})}{\mathcal{H}}(z)\sigma dzd\bar{z}$.
 \end{center}

 Another consequence of formula $(3)$ is the Bochner identity (see \cite {SY}), as now $log{\mathcal{H}}$ 
 is well defined:
 \begin{eqnarray}
 \Delta log{\mathcal{H}}= 2{\mathcal{H}}-2{\mathcal{L}} -2.
 \end{eqnarray} 
 The minimal principle implies ${\mathcal{H}}(z) \ge 1$ for all $z \in \Sigma$, and we find
 \begin{center}
 $\int_{\Sigma}{\mathcal{L}}dA \le \int_{\Sigma}{\mathcal{H}}{\mathcal{L}}dA = 
 \int_{\Sigma}|\mu|^{2}dA = ||\mu||_{WP}= 1$. 
 \end{center}

It is not hard to see that we can actually bound the total energy of this {\hm} $w$ from above, 
in terms of the genus. More precisely, we recall that $w$ is a diffeomorphism, so
\begin{center}
$ \int_{\Sigma}({\mathcal{H}}-{\mathcal{L}})dA =  \int_{\Sigma}JdA = Area(w(\Sigma)) = 4\pi (g-1)$,
\end{center}
therefore the total energy satisfies
\begin{eqnarray*}
E(w) & = & \int_{\Sigma}e(w)dA =  \int_{\Sigma}({\mathcal{H}}+{\mathcal{L}})dA \\
& = & \int_{\Sigma}({\mathcal{H}}-{\mathcal{L}})dA + 2\int_{\Sigma}{\mathcal{L}}dA \\
&\le& Area(\Sigma,\sigma) + 2 \\
&= & 4\pi(g-1) +2.
\end{eqnarray*}
Therefore
 \begin{eqnarray*}
\int_{B_{z_{0}}({\frac{r_{0}}{2}})}{\mathcal{H}}(z)\sigma dzd\bar{z} 
& < & \int_{B_{z_{0}}({\frac{r_{0}}{2}})}({\mathcal{H}}(z) + {\mathcal{L}}(z))\sigma dzd\bar{z} \\
& = & E(w) - \int_{\Sigma \backslash B_{z_{0}}({\frac{r_{0}}{2}})}e(z) \sigma dzd\bar{z} \\
& \le & 4\pi(g-1) +2 - (4\pi(g-1)-A_{1}({B_{z_{0}}}({\frac{r_{0}}{2}}))) \\
&=& A_{1}({B_{z_{0}}}({\frac{r_{0}}{2}})) + 2.
\end{eqnarray*}
where $A_{1}({B_{z_{0}}}({\frac{r_{0}}{2}}))$ is the hyperbolic area of the ball 
$B_{z_{0}}({\frac{r_{0}}{2}})$.

We set $h_0 = C(r_0)(A_{1}({B_{z_{0}}}({\frac{r_{0}}{2}})) + 2)$, then $h_0$ 
is independent of the genus $g$, as it is obtained from a local estimate in a 
geodesic ball. Therefore, 
 \begin{center}
 $|\mu(z)| < {\mathcal{H}}(z) \le {\mathcal{H}}(z_0) < h_0$,
 \end{center}
for all $z \in \Sigma$.
\end{proof}
\begin{rem}
The assumption of the surface lying in the thick part of the {\ms} is essential to this 
argument, since we used an estimate in an embedded geodesic ball.
\end{rem}
\noindent
As an application of this estimate, we can prove theorem 1.1 easily. 

\begin{proof}(of theorem 1.1) Recall that the operator $D = -2(\Delta -2)^{-1}$ is self-adjoint, and $||\mu||_{WP} = 1$:
\begin{center}
$|K_h| = 2\int_{\Sigma}D(|\mu|^2)|\mu|^2 dA < 2h_0^2 \int_{\Sigma}D(|\mu|^2)dA = 2h_0^2$.
\end{center}
\end{proof}.

We now shift our attention to general {\WP} {\Sc}s and to prove theorem 1.2.
\begin{proof}(of theorem 1.2)
We recall from the previous section the Riemannian curvature tensor of the {\WP} metric 
is given by (\cite {Tr}, \cite{Wp3}):
\begin{center} 
$R_{\alpha\bar{\beta}\gamma\bar{\delta}} = 
\int_{\Sigma}D(\mu_{\alpha}\bar\mu_{\beta})\mu_{\gamma}\bar\mu_{\delta}dA + 
\int_{\Sigma}D(\mu_{\alpha}\bar\mu_{\delta})\mu_{\gamma}\bar\mu_{\beta}dA$,
\end{center}
where $\mu$'s in the formula are harmonic {\Bd}s, and $D$ again is the operator 
$-2(\Delta -2)^{-1}$.

To calculate the {\Sc}, we choose two arbitrary orthonormal harmonic {\Bd}s 
$\mu_0$ and $\mu_1$. In other words, we have 
\begin{center}
$\int_{\Sigma}|\mu_0|^2 dA = \int_{\Sigma}|\mu_1|^2 dA = 1$ and 
$\int_{\Sigma}\mu_0 {\bar{\mu}}_1 dA = 0$.
\end{center}
Then the Gaussian curvature of the plane spanned by $\mu_0$ and $\mu_1$ is 
(\cite {Wp3})
\begin{eqnarray*}
K(\mu_0,\mu_1) &=& {\frac{1}{4}}(R_{0\bar{1}0\bar{1}} - R_{0\bar{1}1\bar{0}} - 
R_{1\bar{0}0\bar{1}}+R_{1\bar{0}1\bar{0}}) \nonumber \\
&=& Re(\int_{\Sigma}D(\mu_{0}\bar\mu_{1})\mu_{0}\bar\mu_{1}dA) - 
{\frac{1}{2}}Re(\int_{\Sigma}D(\mu_{0}\bar\mu_{1})\mu_{1}\bar\mu_{0}dA) \nonumber \\
&-& {\frac{1}{2}}\int_{\Sigma}D(|\mu_{1}|^2)|\mu_{0}|^2dA. 
\end{eqnarray*}
From (\cite {Wp3}, lemma 4.3) and H\"{o}lder inequality, we have
\begin{eqnarray*}
|Re(\int_{\Sigma}D(\mu_{0}\bar\mu_{1})\mu_{0}\bar\mu_{1}dA)| &\le& 
\int_{\Sigma}|D(\mu_{0}\bar\mu_{1})||\mu_{0}\bar\mu_{1}|dA \nonumber \\
&\le& \int_{\Sigma}\sqrt{D(|\mu_{0}|^{2})}\sqrt{D(|\mu_{0}|^{2})}|\mu_{0}\bar\mu_{1}|dA \nonumber \\
&\le&  \int_{\Sigma}D(|\mu_{0}|^2)|\mu_{1}|^2dA = \int_{\Sigma}D(|\mu_{1}|^2)|\mu_{0}|^2dA, 
\end{eqnarray*}
and similarly
\begin{center}
$|\int_{\Sigma}D(\mu_{0}\bar\mu_{1})\mu_{1}\bar\mu_{0}dA| \le \int_{\Sigma}D(|\mu_{1}|^2)|\mu_{0}|^2dA$. 
\end{center}
Therefore 
\begin{center}
$|K(\mu_0,\mu_1)| \le 2\int_{\Sigma}D(|\mu_{1}|^2)|\mu_{0}|^2dA$.  
\end{center}
We apply theorem 3.3 to find that there is a $h_0 > 0$, independent of $g$, such that $|\mu_{0}| < h_0$. 
So $|K| < 2h_0^2 \int_{\Sigma}D(|\mu_{1}|^2)dA = 2h_0^2$.
\end{proof}

It is now straightforward to see that theorem 1.3 holds since all {\Sc}s are negative.

\end{document}